\documentclass{amsart}
\usepackage{amssymb}
\def\version{Version W02-v3c     last changed 31 March 2004 by SN}
\date{\version}
\usepackage{a4}
\usepackage[notref,notcite]{showkeys}
\begin{document}
\newtheorem{theorem}{Theorem}[section]
\newtheorem{lemma}[theorem]{Lemma}
\newtheorem{example}[theorem]{Example}
\def\ffrac#1#2{{\textstyle\frac{#1}{#2}}}
\def\AA{A^1}\def\BB{{\mathcal{B}}}
\makeatletter
 \renewcommand{\theequation}{%
 \thesection.\alph{equation}}
 \@addtoreset{equation}{section}
 \makeatother
\title[Curvature homogeneous manifolds]
{Complete curvature homogeneous pseudo-Riemannian manifolds}
\author{P. Gilkey and S. Nik\v cevi\'c}
\begin{address}{PG: Mathematics Department, University of Oregon,
Eugene Or 97403 USA.\newline Email: {\it gilkey@darkwing.uoregon.edu}}
\end{address}
\begin{address}{SN: Mathematical Institute, Sanu,
Knez Mihailova 35, p.p. 367,
11001 Belgrade,
Serbia and Montenegro.
\newline Email: {\it stanan@mi.sanu.ac.yu}}\end{address}
\begin{abstract} We exhibit 3 families of complete curvature homogeneous
pseudo-Riemannian manifolds which are modeled on irreducible symmetric
spaces and which are not locally homogeneous. All of the manifolds have nilpotent Jacobi operators; some of the manifolds
are, in addition, Jordan Osserman and Jordan Ivanov-Petrova.\end{abstract}
\maketitle

\section{Introduction} Consider a triple $\mathcal{U}:=(V,g,A)$ where $g$
is a non-degenerate inner product of signature $(p,q)$ on an
$m$-dimensional real vector space $V$ with $m=p+q$ and where
$A\in\otimes^4V^*$ is an algebraic curvature tensor -- i.e. a
$4$-tensor satisfying the usual symmetries of the Riemann
curvature tensor:
\begin{eqnarray*}
&&A(x,y,z,w)=-A(y,x,z,w)=A(z,w,x,y)\ \ \text{and}\\
&&A(x,y,z,w)+A(y,z,x,w)+A(z,x,y,w)=0\,.
\end{eqnarray*}

We also consider a pair $\mathcal{M}:=(M,g_M)$ where $g_M$ is a pseudo-Riemannian metric
 of signature $(p,q)$ on a manifold $M$ of dimension $m=p+q$. One says $\mathcal{M}$ is {\it Riemannian} if $p=0$ and {\it
Lorentzian} if $p=1$. Let $R_M$ be the associated Riemann curvature tensor. We say that
$\mathcal{U}$ is a {\it $0$-model for}
$\mathcal{M}$ if for every point $P\in M$, there exists an isomorphism $\Phi_P:T_PM\rightarrow V$ so that
$$\Phi_P^*g=g_M|_{T_PM}\ \ \text{and}\ \ \Phi_P^*A=R_M|_{T_PM}\,.$$
One says that $\mathcal{M}$ is {\it curvature homogeneous}
if
$\mathcal{M}$ admits a $0$-model, in other words, the metric and curvature tensor
``look the same at each point''. If $\mathcal{N}:=(N,g_N)$ is a homogeneous space, we say that $\mathcal{M}$ is
{\it modelled on}
$\mathcal{N}$ if $(T_QN,h_N|_{T_QN},R_N|_{T_QN})$ is a $0$-model for $(M,g_M,R_M)$; the precise $Q\in N$
being immaterial since $\mathcal{N}$ is assumed to be homogeneous. We refer to \cite{KTV92,T88,Va91} for further details.

We say that $\AA\in\otimes^5V^*$ is an {\it algebraic covariant derivative
curvature tensor} if
$\AA$ has the curvature symmetries of the covariant derivative of the Riemann curvature tensor, i.e. we have the relations:
\begin{eqnarray*}
&&\AA(x,y,z,w;v)=\AA(z,w,x,y;v)=-\AA(y,x,z,w;v),\nonumber\\
&&\AA(x,y,z,w;v)+\AA(y,z,x,w;v)+\AA(z,x,y,w;v)=0,\\
&&\AA(x,y,z,w;v)+\AA(x,y,w,v;z)+\AA(x,y,v,z;w)=0.\nonumber
\end{eqnarray*}
We say that a quadruple $\mathcal{U}^1:=(V,g,A,\AA)$ is a {\it 1-model} for $\mathcal{M}$ if for every point $P\in M$, there
exists an isomorphism
$\Phi_P:T_PM\rightarrow V$ so that
$$
\Phi_P^*g=g_M|_{T_PM},\ \ \Phi_P^*A=R_M|_{T_PM},\ \ \text{and}\ \
\Phi_P^*\AA=\nabla R_M|_{T_PM}\,.
$$
In this setting, $\mathcal{M}$ is said to be {\it $1$-curvature homogeneous}.
The notion of $k$-curvature homogeneous for $k\ge2$ is defined similarly. These notions were first introduced by Singer who
showed:
\begin{theorem}\label{thm-1.1}
{\bf (Singer \cite{S60})} There exists an universal bound $k_m$
such that a Riemannian manifold $\mathcal{M}$ of dimension $m$ is
locally homogeneous if and only if $\mathcal{M}$ is
$(k_m+1)$-curvature homogeneous. Furthermore, $k_m$ is smaller
than $m(m-1)/2$.
\end{theorem}

One has the following important results in the context of models based on the curvature tensors of irreducible symmetric spaces
in the Riemannian and Lorentzian setting:

\begin{theorem}\label{thm-1.2}
\ \begin{enumerate}
\item {\bf (Tricerri and Vanhecke \cite{TV})}
A Riemannian curvature homogeneous manifold modelled on an
irreducible symmetric space is locally symmetric.
\item {\bf (Cahen, Leroy, Parker, Tricerri, and Vanhecke \cite{CLPT90})} A Lorentzian curvature homogeneous manifold
modelled on an irreducible symmetric space has constant sectional curvature.\end{enumerate}\end{theorem}

The proof of Theorem \ref{thm-1.2} (1) uses properties of the scalar curvature invariants of
Riemannian manifolds which do not hold for indefinite metrics; the proof of Theorem \ref{thm-1.2} (2) uses the
remark, which is based on M. Berger's classification, that any irreducible Lorentzian symmetric space of dimension greater
than or equal to 3 has constant sectional curvature. In this brief
note, we will present several examples illustrating that Theorem \ref{thm-1.2} fails in the
higher signature setting by constructing complete curvature homogeneous pseudo-Riemannian manifolds which are modelled on
irreducible symmetric spaces and which are not locally homogeneous.

Throughout this paper, we will be introducing metrics, curvature tensors, and covariant derivative curvature tensors.
In the interests of brevity, we shall often only give the non-zero components of these tensors up to the usual $\mathbb{Z}_2$
symmetries.

\subsection{Signature $(p,p)$}\label{sect-1.1} There are curvature homogeneous pseudo-Riemannian
manifolds of balanced (or neutral) signature $(p,p)$ which are complete but not locally homogeneous (and
hence not locally symmetric), but which nevertheless are modeled on a complete irreducible symmetric
space.

Let $p\ge3$. Let $(\vec x,\vec y)$ for
$\vec x=(x_1,...,x_p)$ and
$\vec y=(y_1,...,y_p)$ be coordinates on $\mathbb{R}^{2p}$. Let
$f=f(\vec x)$ be a smooth function on $\mathbb{R}^p$. Let
$\mathcal{M}_{1,p,f}:=(\mathbb{R}^{2p},g_{1,p,f})$ where
$$
g_{1,p,f}(\partial_i^x,\partial_j^x)=\partial_i^xf\cdot\partial_j^xf\ \ \text{and}\ \
g_{1,p,f}(\partial_i^x,\partial_i^y)=1,
$$
the other components being zero.
Let $\mathbb{R}^{2p}=\operatorname{Span}\{X_1,...,X_p,Y_1,...,Y_p\}$. Let
$\mathcal{U}_{1,p}:=(\mathbb{R}^{2p},g_{1,p},A_{1,p})$ where
$$
g_{1,p}(X_i,Y_i)=1\ \ \text{and}\ \
A_{1,p}(X_i,X_j,X_k,X_l)=\delta_{il}\delta_{jk} -\delta_{ik}\delta_{jl}\,.
$$
The metric $g_{1,p,f}$ and the inner product $g_{1,p}$ have signature $(p,p)$.

Let
$H_f=(H_{f,ij})$, where $H_{f,ij}:=(\partial_i^x\partial_j^xf)$, be the Hessian matrix of second partial derivatives. Assume
$H_f>0$. Let $H_f^{ij}$ be the inverse matrix. Let $R_{1,p,f}$ be the curvature tensor of the metric $g_{1,p,f}$ and let
$\nabla R_{1,p,f}$ be the associated covariant derivative. Set
\begin{eqnarray*}
&&\alpha_1:=\textstyle\sum_{a,b,c,d,e,s,t,u,v,w}H_f^{as}H_f^{bt}H_f^{cu}H_f^{dv}H_f^{ew}
\nabla
R_{1,p,f}(\partial_{a}^x,\partial_{b}^x,\partial_{c}^x,\partial_{d}^x;\partial_{e}^x)\\
&&\qquad\qquad\qquad\qquad\qquad\cdot
\nabla R_{1,p,f}(\partial_{s}^x,\partial_{t}^x,\partial_{u}^x,\partial_{v}^x;\partial_{w}^x)
\,.\end{eqnarray*}

\begin{theorem}\label{thm-1.3}
 Let $p\ge2$ and let $H_f>0$. Then:
\begin{enumerate}
\item All geodesics in $\mathcal{M}_{1,p,f}$ extend for infinite time.
\item If $P\in\mathbb{R}^{2p}$, then $\exp_P:T_P\mathbb{R}^{2p}\rightarrow\mathbb{R}^{2p}$ is
a diffeomorphism.
\item The non-zero components of $R_{1,p,f}$ and of $\nabla R_{1,p,f}$ are given by:
\begin{enumerate}\item $R_{1,p,f}(\partial_i^x,\partial_j^x,\partial_k^x,\partial_l^x)=
H_{f,il}H_{f,jk}-H_{f,ik}H_{f,jl}$,
\item $\nabla R_{1,p,f}(\partial_i^x,\partial_j^x,\partial_k^x,\partial_l^x;\partial_n^x)=
\partial_n^xR_{1,p,f}(\partial_i^x,\partial_j^x,\partial_k^x,\partial_l^x)$.
\end{enumerate}
\item $\mathcal{U}_{1,p}$ is an irreducible $0$-model for $\mathcal{M}_{1,p,f}$.
\item If $f=x_1^2+...+x_p^2$, then $\mathcal{M}_{1,p,f}$ is
an irreducible symmetric space.
\item If $p\ge3$ and if $\alpha_1$ is not constant, $\mathcal{M}_{1,p,f}$ is not curvature $1$-homogeneous.
\end{enumerate}
\end{theorem}

The pseudo-Riemannian manifold $\mathcal{M}_{1,p,f}$ can be realized as a hypersurface in a flat space of signature $(p,p+1)$;
the Hessian $H_f$ then gives the second fundamental form. We refer to \cite{refDG,GIZ02,GIZ03} for further details concerning
this family of manifolds.

\subsection{Signature $(2s,s)$}\label{sect-1.2} For $s\ge2$, let $(\vec  u,\vec t,\vec v)$ give coordinates on
$\mathbb{R}^{3s}$ where we have
$\vec u:=(u_1,...,u_s)$,
$\vec t:=(t_1,...,t_s)$, and
$\vec v:=(v_1,...,v_s)$.
Let $f_i\in C^\infty(\mathbb{R})$ be smooth functions. Set
$$F(\vec u):=f_1(u_1)+...+f_s(u_s)\in C^\infty(\mathbb{R}^s)\ \ \text{and}\ \ |u|^2=u_1^2+...+u_s^2\,.$$
Let $\mathcal{M}_{2,s,F}:=(\mathbb{R}^{3s},g_{2,s,F})$ where
\begin{eqnarray*}
&&g_{2,s,F}(\partial_i^u,\partial_j^u)=-2\{F(\vec u)+\textstyle\sum_{1\le k\le s}u_kt_k\}\delta_{ij},\\
&&g_{2,s,F}(\partial_i^u,\partial_j^v):=\delta_{ij},\ \ \text{and}\ \
 g_{2,s,F}(\partial_i^t,\partial_j^t):=-\delta_{ij}\,.\end{eqnarray*}
Manifolds of this type were first introduced in \cite{GN03x}; see also \cite{GN03,GNV03}.

Let $\mathbb{R}^{3s}=\operatorname{Span}\{U_1,...,U_s,T_1,...,T_s,V_1,...,V_s\}$. Let
$\mathcal{U}_{2,s}:=(\mathbb{R}^{3s},g_{2,s},A_{2,s})$ for
\begin{equation}\label{eqn-1.a}
\begin{array}{l}
g_{2,s}(U_i,V_j):=\delta_{ij},\ \ g_{2,s}(T_i,T_j):=-\delta_{ij},\
\ \text{and}\\
A_{2,s}(U_i,U_j,U_k,T_l):=\delta_{il}\delta_{jk}-\delta_{ik}\delta_{jl}\,.
\vphantom{\vrule height 12pt}\end{array}\end{equation} The metric
$g_{2,s,F}$ and the inner product $g_{2,s}$ have signature
$(2s,s)$. Set $$\alpha_2:=\textstyle\sum_{i,j,k,l,n}\{\nabla
R_{2,s,F}(\partial_i^u,\partial_j^u,\partial_k^u,\partial_l^u;\partial_n^u)\}^2\,.$$

\begin{theorem}\label{thm-1.4}
Let $s\ge2$. Then:\begin{enumerate}
\item All geodesics in $\mathcal{M}_{2,s,F}$ extend for infinite time.
\item If $P\in\mathbb{R}^{3s}$, then $\exp_P:T_P\mathbb{R}^{3s}\rightarrow\mathbb{R}^{3s}$ is a diffeomorphism.
\item The non-zero components of $R_{2,s,F}$ and of $\nabla R_{2,s,F}$ are given by:
\begin{enumerate}\item $R_{2,s,F}(\partial_i^u,\partial_j^u,\partial_j^u,\partial_i^u)
 =(\partial_i^u)^2f_i+(\partial_j^u)^2f_j+|u|^2$,
\item $R_{2,s,F}(\partial_i^u,\partial_j^u,\partial_j^u,\partial_i^t)=1$,
\item $\nabla R_{2,s,F}(\partial_i^u,\partial_j^u,\partial_j^u,\partial_i^u;\partial_i^u)
 =(\partial_i^u)^3f_i+4u_i$.
\end{enumerate}
\item $\mathcal{U}_{2,s}$ is an irreducible $0$-model for $\mathcal{M}_{2,s,F}$.
\item If $F=-\frac16u_1^4-...-\frac16u_s^4$, then
$\mathcal{M}_{2,s,F}$ is an irreducible symmetric space.
\item If $s\ge3$ and if $\alpha_2$ is not constant, $\mathcal{M}_{2,s,F}$ is not curvature $1$-homogeneous.
\end{enumerate}
\end{theorem}

Assertion (6) in
Theorem \ref{thm-1.4} was discussed previously in \cite{GN03}; C. Dunn pointed out that the argument given there contained a
mistake. In this paper, we shall give a slightly different argument which avoids that mistake.

\subsection{Manifolds which are $1$-curvature homogeneous}\label{sect-1.3} The previous two families of examples were
curvature homogeneous but not $1$-curvature homogeneous for generic members of the families. Let
$r\ge2$. Introduce coordinates
$(\vec u,\vec v,x,y)$ on
$\mathbb{R}^{2r+2}$ where we have
$\vec u=(u_1,...,u_r)$ and
$\vec v=(v_1,...,v_r)$. Let $\psi\in C^\infty(\mathbb{R})$. Let
$\mathcal{M}_{3,r,\psi}:=(\mathbb{R}^{2r+2},g_{3,r,\psi})$ where
\begin{eqnarray*}
&&g_{3,r,\psi}(\partial_x,\partial_y)=1,\ \
 g_{3,r,\psi}(\partial_{u_i},\partial_{v_j})=\delta_{ij},\ \ \text{and}\\
&&g_{3,r,\psi}(\partial_x,\partial_x)=-2u_1v_2-...-2u_{r-1}v_r-2\psi(u_r)\,.
\end{eqnarray*}
These manifolds are closely related to examples of
Fiedler et al \cite{refFeGi03}.

Let $\mathbb{R}^{2r+2}=\operatorname{Span}\{U_1,...,U_r,V_1,...,V_r,X,Y\}$.
Let $\mathcal{U}_{3,r}:=(\mathbb{R}^{2r+2},g_{3,r},A_{3,r})$ for
\begin{eqnarray*}
&&\ \ g_{3,r}(X,Y)=1,\ \
g_{3,r}(U_i,V_j)=\delta_{ij},\ \
A_{3,r}(X,U_r,U_r,X)=1,\ \ \text{and}\\&&
A_{3,r}(X,U_i,V_{i+1},X)=1\ \ \text{for}\ \ 1\le i\le r-1\,.
\end{eqnarray*}
The metric $g_{3,r,\psi}$ and the inner product $g_{3,r}$ have signature $(r+1,r+1)$. We also define a
$1$-model space $\mathcal{U}_{3,r}^{1}:=(\mathbb{R}^{2r+2},g_{3,r},A_{3,r},\AA_{3,r})$ where
$$
\AA_{3,r}(X,U_r,U_r,X;U_r)=1\,.
$$

\begin{theorem}\label{thm-1.5}
Let $r\ge2$. Assume that $\psi^{\prime\prime}>0$. Then:
\begin{enumerate}
\item All geodesics in $\mathcal{M}_{3,r,\psi}$ extend for infinite time.
\item $\exp_P:T_P\mathbb{R}^{2r+2}\rightarrow\mathbb{R}^{2r+2}$ is a diffeomorphism for all $P$ in
$\mathbb{R}^{2r+2}$.
\item The non-zero components of $R_{3,r,\psi}$ and of $\nabla R_{3,r,\psi}$ are given by:
\begin{enumerate}
\item $R_{3,r,\psi}(\partial_x,\partial_{u_r},\partial_{u_r},\partial_x)=\psi^{\prime\prime}(u_r)$,
\item $R_{3,r,\psi}(\partial_x,\partial_{u_i},\partial_{v_{i+1}},\partial_x)=1$ for $1\le i\le r-1$,
\item $\nabla R_{3,r,\psi}(\partial_x,\partial_{u_r},\partial_{u_r},\partial_x;\partial_{u_r})=\psi^{\prime\prime\prime}(u_r)$.
\end{enumerate}
\item $\mathcal{U}_{3,r}$ is an irreducible $0$-model for $\mathcal{M}_{3,r,\psi}$.
\item If $\psi(u_r)=u_r^2$, then $\mathcal{M}_{3,r,\psi}$ is an irreducible symmetric space.
\item If $\psi^{\prime\prime\prime}>0$, then:\begin{enumerate}\item $\mathcal{U}_{3,r}^1$ is a
$1$-model for $\mathcal{M}_{3,r,\psi}$,
\item $\mathcal{M}_{3,r,\psi}$ is not $2$-curvature homogeneous.
\end{enumerate}\end{enumerate}
\end{theorem}

Theorem \ref{thm-1.3} (6) (resp. Theorem \ref{thm-1.4} (6)) requires that $p\ge3$ (resp. $s\ge3$). This result is sharp; for
suitably chosen $f$ (resp. $F$), $\mathcal{M}_{1,2,f}$ (resp. $\mathcal{M}_{2,2,F}$) is curvature $1$-homogeneous
but not curvature
$2$-homogeneous; we omit details in the interests of brevity.

\subsection{Osserman manifolds}\label{sect-1.4} If $x$ is a tangent vector at a point $P\in M$, then the
{\it Jacobi operator} $J_M(x)$ is characterized by the identity
$$g_M(J_M(x)y,z)=R_M(y,x,x,z)\,.$$ If $\rho_M$ is the associated
Ricci tensor, then $\rho_M(x,x)=\operatorname{Tr}(J_M(x))$. One
says that $\mathcal{M}$ is {\it spacelike} (resp. {\it timelike})
{\it Osserman} if the eigenvalues of the Jacobi operator are
constant on the pseudo-sphere bundle $S^+(\mathcal{M})$ of
spacelike (resp. $S^-(\mathcal{M})$ of timelike) unit vectors. One
says that $\mathcal{M}$ is {\it spacelike} (resp. {\it timelike})
{\it Jordan Osserman} if the Jordan normal form of the Jacobi
operator is constant on $S^+(\mathcal{M})$ (resp.
$S^-(\mathcal{M})$). We shall say that $\mathcal{M}$ is {\it
Osserman nilpotent of order n} if $J_M(x)^n=0$ for every $x\in TM$
and if there exists a point $P_0\in M$ and a tangent vector
$x_0\in T_{P_0}M$ so that $J_M(x_0)^{n-1}\ne0$. Such manifolds are
necessarily Osserman since $0$ is the only eigenvalue of $J_M$.
And consequently such manifolds are Ricci flat since
$\rho(x,x)=\operatorname{Tr}(J(x))$. We refer to
\cite{GKV02,refGi02} for further details concerning Osserman
manifolds.

\begin{theorem}\label{thm-1.6}
\ \begin{enumerate}
\item Let $p\ge2$. If $H_f>0$, then
$\mathcal{M}_{1,p,f}$ is spacelike and timelike Jordan Osserman.
\item Let $s\ge2$. Then $\mathcal{M}_{2,s,F}$ is spacelike Jordan Osserman. However $\mathcal{M}_{2,s,F}$ is not timelike
Jordan Osserman.
\item Let $r\ge2$. If $\psi^{\prime\prime}>0$, then $\mathcal{M}_{3,r,\psi}$ is $2r$-Osserman
nilpotent.
\end{enumerate}\end{theorem}

The three families $\mathcal{M}_{i,k}$ first arose in the study of Osserman manifolds. We refer to
\cite{BBGR97,BV03,GKV02,GVV98} for other examples of non-homogeneous Osserman manifolds.

\subsection{Ivanov-Petrova manifolds} Let $\{e_1,e_2\}$ be an oriented orthonormal basis for an oriented spacelike (resp.
timelike) $2$-plane $\pi$. The {\it skew-symmetric curvature operator} $\mathcal{R}_M(\pi)$ is characterized by the identity
$$g_M(\mathcal{R}_M(\pi)y,z)=R_M(e_1,e_2,y,z)\,.$$
This operator is independent of the particular oriented orthonormal basis chosen for $\pi$. One says that
$\mathcal{M}$ is {\it spacelike} (resp. {\it timelike}) {\it Jordan Ivanov-Petrova} if the Jordan normal form of
$\mathcal{R}_M$ is constant on the Grassmannian of oriented spacelike (resp. timelike)
$2$-planes; one lets the {\it Rank} be the common value of $\operatorname{rank}(\mathcal{R}_M(\pi))$ in this setting.

\begin{theorem}\label{thm-1.7}\
\begin{enumerate}
\item Let $p\ge2$. If $H_f>0$, then $\mathcal{M}_{1,p,f}$ is
spacelike and timelike Jordan Ivanov-Petrova of rank $2$.
\item Let $s\ge2$. Then $\mathcal{M}_{2,s,F}$ is spacelike
Jordan Ivanov-Petrova of rank $4$; $\mathcal{M}_{2,s,F}$ is not timelike Jordan Ivanov-Petrova.
\end{enumerate}\end{theorem}

\subsection{The geodesic involution}\label{sect-1.6}
The following observation is a
special case which follows from work of E. Cartan; we present it for the sake of
completeness in light of the examples given above.

\begin{theorem}\label{thm-1.8}
 Let $\mathcal{M}$ be a pseudo-Riemannian manifold of signature
$(p,q)$. Suppose that
$\nabla R_M=0$ and that
$\exp_P:T_PM\rightarrow M$ is a diffeomorphism for every $P\in M$.
Then the geodesic symmetry $\mathcal{S}_P:Q\rightarrow
\exp_P\{-\exp_P^{-1}Q\}$ is an isometry. Furthermore, $\mathcal{M}$ is a homogeneous space.
\end{theorem}

Here is a brief guide to this paper. In Section \ref{sect-2}, we prove Assertions (1)-(3) of Theorems
\ref{thm-1.3}-\ref{thm-1.5}. In Section \ref{sect-3}, we show $\mathcal{U}_{i,k}$ is a $0$-model for
$\mathcal{M}_{i,k}$ and in Section \ref{sect-4}, we show these models are irreducible. This establishes
Assertion (4) of Theorems \ref{thm-1.3}-\ref{thm-1.5}. Assertion (5) of these three Theorems then follows as a scholium to the
previous assertions. In Section \ref{sect-5}, we establish Assertion (6) of Theorems \ref{thm-1.3}-\ref{thm-1.5}. We refer to
\cite{GIZ03} for the proof of Assertion (1) of Theorems \ref{thm-1.6} and \ref{thm-1.7} and to \cite{GN03,GNV03} for the
proof of Assertion (2) of Theorems \ref{thm-1.6} and \ref{thm-1.7}. Assertion (3) of Theorem \ref{thm-1.6} is proved in Section
\ref{sect-4}. In Section \ref{sect-6}, we complete our discussion by proving Theorem
\ref{thm-1.8}.

It is a pleasant task to thank
Professors E. Garc\'{\i}a--R\'{\i}o, O. Kowalski, and L. Vanhecke for useful conversations on this subject.

\section{Complete manifolds}\label{sect-2}

We shall need the following technical fact.
\begin{lemma}\label{lem-2.1}
Let $(z_1,...,z_n)$ be coordinates on $\mathbb{R}^n$. Let $g$ be a
pseudo-Riemannian metric on $\mathbb{R}^n$ so that
$\nabla_{\partial_a^z}\partial_b^z=\sum_{a,b<c}\Gamma_{ab}{}^c(z_1,...,z_{c-1})\partial_c^z$.
Then:\begin{enumerate}
\item $(\mathbb{R}^n,g)$ is a complete pseudo-Riemannian manifold.
\item $\exp_P:T_P\mathbb{R}^n\rightarrow\mathbb{R}^n$ is a diffeomorphism for all $P$ in $\mathbb{R}^n$.
\end{enumerate}
\end{lemma}

\begin{proof} Let $\gamma(t)=(z_1(t),...,z_n(t))$ be a curve in $\mathbb{R}^n$; $\gamma$ is a
geodesic if and only
\begin{eqnarray*}
&&\ddot z_1(t)=0,\ \ \text{and for}\ \ c>1\\
&&\ddot z_c(t)+\textstyle\sum_{a,b<c}\dot z_a(t)\dot z_b(t)\Gamma_{ab}{}^c(z_1,...,z_{c-1})(t)=0\,.
\end{eqnarray*}
We solve this system of equations recursively. Let $\gamma(t;\vec z^{\phantom{.}0},\vec z^{\phantom{.}1})$ be defined by
\begin{eqnarray*}
&&z_1(t)=z_1^0+z_1^1t,\ \ \text{and for}\ \ c>1\\
&&z_c(t)=z_c^0+z_c^1t-\textstyle\int_0^t\int_0^s\textstyle\sum_{a,b<c}\dot z_a(r)\dot
z_b(r)\Gamma_{ab}{}^c(z_1,...,z_{c-1})(r)drds\,.
\end{eqnarray*}
Then $\gamma(0;\vec z^{\phantom{.}0},\vec z^{\phantom{.}1})(0)=\vec z^{\phantom{.}0}$ while $\dot\gamma(0;\vec
z^{\phantom{.}0},\vec z^{\phantom{.}1})(0)=\vec z^{\phantom{.}1}$. Thus every geodesic arises in this way so all geodesics
extend for infinite time. Furthermore, given $P,Q\in\mathbb{R}^n$, there is a unique geodesic $\gamma=\gamma_{P,Q}$ so that
$\gamma(0)=P$ and
$\gamma(1)=Q$ where
\begin{eqnarray*}
&&z_1^0=P_1,\ \ z_1^1=Q_1-P_1,\ \ \text{and for}\ \ c>1\\
&&z_c^0=P_c,\ \ z_c^1=Q_c-P_c+\textstyle\int_0^1\int_0^s\textstyle\sum_{a,b<c}\dot
z_a(r)\dot z_b(r)\Gamma_{ab}{}^c(z_1,...,z_{c-1})(r)drds\,.
\end{eqnarray*}
This shows that $\exp_P$ is a diffeomorphism from $T_P\mathbb{R}^n$ to $\mathbb{R}^n$.
\end{proof}

\begin{proof}[Proof of Theorem \ref{thm-1.3} (1-3)] Adopt the notation of Section \ref{sect-1.1}. Let
$$g_{ij}(x)=g_{1,p,f}(\partial_i^x,\partial_j^x):=\partial_i^xf\cdot\partial_j^xf,\ \ \text{and}\ \
\Gamma_{ijk}(\vec x):=\ffrac12\{\partial_i^xg_{jk}+\partial_j^xg_{ik}-
\partial_k^xg_{ij}\}\,.$$
The non-zero Christoffel symbols are
\begin{equation}\label{eqn-2.a}\begin{array}{lll}
g_{1,p,f}(\nabla_{\partial_i^x}\partial_j^x,\partial_k^x)=\Gamma_{ijk}(\vec x)\ \ \text{and}\ \
\nabla_{\partial_i^x}\partial_j^x=\textstyle\sum_{1\le k\le p}\Gamma_{ijk}(\vec x)\partial_k^y\,.
\end{array}\end{equation}
We verify that the hypothesis of Lemma \ref{lem-2.1} is satisfied and thereby prove Assertions (1) and (2) by setting:
$$z_1=x_1,\ ...,\ z_p=x_p,\ z_{p+1}=y_1,\ ...,\ z_{2p}=y_p\,.$$
Furthermore, by Equation (\ref{eqn-2.a}),
\begin{eqnarray*}R_{1,p,f}(\partial_i^x,\partial_j^x,\partial_k^x,\partial_l^x)&=&
-\ffrac12(\partial_j^x\partial_k^xg_{il}+\partial_i^x\partial_l^xg_{jk}
 -\partial_j^x\partial_l^xg_{ik}-\partial_i^x\partial_k^xg_{jl})\\
&=&H_{f,il}H_{f,jk}-H_{f,ik}H_{f,jl}
\end{eqnarray*}
while $R_{1,p,f}(\cdot,\cdot,\cdot,\cdot)=0$ if any of the entries is $\partial_i^y$. Assertion (3a) now follows. Furthermore, by
Equation (\ref{eqn-2.a}),
$$\nabla R_{1,p,f}(\partial_i^x,\partial_j^x,\partial_k^x,\partial_l^x;\partial_n^x)=
\partial_n^xR_{1,p,f}(\partial_i^x,\partial_j^x,\partial_k^x,\partial_l^x)$$
while $\nabla R_{1,p,f}(\cdot,\cdot,\cdot,\cdot;\cdot)=0$ if any of the entries is $\partial_i^y$.
This proves Assertion (3b) of Theorem \ref{thm-1.3}.
\end{proof}

\begin{proof}[Proof of Theorem \ref{thm-1.4} (1)-(3)] Adopt the notation of Section \ref{sect-1.2}. Let
$i\ne j$ and let $g=g_{2,s,F}$. The non-zero Christoffel symbols of the second kind are given by:
$$\begin{array}{ll}
g(\nabla_{\partial_i^u}\partial_i^u,\partial_i^u)=-\partial_i^uf_i-t_i,& \\
g(\nabla_{\partial_i^u}\partial_i^u,\partial_j^u)=\partial_j^uf_j+t_j,&
g(\nabla_{\partial_i^u}\partial_j^u,\partial_i^u)=
g(\nabla_{\partial_j^u}\partial_i^u,\partial_i^u)=-\partial_j^uf_j-t_j,\vphantom{\vrule height 11pt}\\
g(\nabla_{\partial_i^u}\partial_i^u,\partial_i^t)=u_i,&
g(\nabla_{\partial_i^u}\partial_i^t,\partial_i^u)=g(\nabla_{\partial_i^t}\partial_i^u,\partial_i^u)=-u_i,\\
g(\nabla_{\partial_i^u}\partial_i^u,\partial_j^t)=u_j,&
g(\nabla_{\partial_i^u}\partial_j^t,\partial_i^u)=
g(\nabla_{\partial_j^t}\partial_i^u,\partial_i^u)=-u_j\,.\vphantom{\vrule height 11pt}
\end{array}$$
We may then raise indices to see the non-zero covariant derivatives are given by:
\begin{eqnarray*}
&&\nabla_{\partial_i^u}\partial_i^u=-(\partial_i^uf_i+t_i)\partial_i^v
 +\textstyle\sum_{k\ne i,1\le k\le s}(\partial_k^uf_k+t_k)\partial_k^v-\textstyle\sum_{1\le k\le s}u_k\partial_k^t,\\
&&\nabla_{\partial_i^u}\partial_j^u=-(\partial_j^uf_j+t_j)\partial_i^v-(\partial_i^uf_i+t_i)\partial_j^v,\\
&&\nabla_{\partial_i^u}\partial_i^t=\nabla_{\partial_i^t}\partial_i^u=-u_i\partial_i^v,\ \ \text{and}\ \
 \nabla_{\partial_i^u}\partial_j^t=\nabla_{\partial_j^t}\partial_i^u=-u_j\partial_i^v\,.
\end{eqnarray*}
We derive Assertions (1) and (2) from Lemma \ref{lem-2.1} by setting:
$$z_1=u_1,\ ...,\ z_s=u_s,\ z_{s+1}=t_1,\ ...,\ z_{2s}=t_p,\ z_{2s+1}=v_1,\ ...,\ z_{3s}=v_s\,.$$

We have $\nabla\partial_i^v=0$. Thus if at least one $z_\mu\in\{\partial_i^v\}$, then
$R_{2,s,F}(z_1,z_2,z_3,z_4)=0$. Similarly, if at least two of the $z_\mu$ belong to $\{\partial_i^t\}$, then
$R_{2,s,F}(z_1,z_2,z_3,z_4)=0$. Finally, as
$\partial_i^u\partial_j^uF=0$ for $i\ne j$,
$R_{2,s,F}(\partial_i^u,\partial_j^u,\partial_k^u,\star)=0$ if the indices $\{i,j,k\}$ are distinct. Furthermore
$$
\nabla_{\partial_i^u}\nabla_{\partial_j^u}\partial_j^u=f_i^{\prime\prime}\partial_i^v-\partial_i^t+\{|u|^2\}\partial_i^v\quad\text{and}\quad
\nabla_{\partial_j^u}\nabla_{\partial_i^u}\partial_j^u=-f_j^{\prime\prime}\partial_i^v\,.
$$
Assertions (3a) and (3b) now follow.

We have similarly that $\nabla R_{2,s,F}(\xi_1,\xi_2,\xi_3,\xi_4;\xi_5)=0$ if at least one of the $\xi_i$ belongs to
$\operatorname{Span}\{T_i,V_i\}$. Furthermore, the only non-zero component of $\nabla R_{2,s,F}$ is given by:
\begin{eqnarray*}
&&\nabla R_{2,s,F}(\partial_i^u,\partial_j^u,\partial_j^u,\partial_i^u;\partial_i^u)\\
&=&\partial_i^uR_{2,s,F}(\partial_i^u,\partial_j^u,\partial_j^u,\partial_i^u)-2
 R_{2,s,F}(\nabla_{\partial_i^u}\partial_i^u,\partial_j^u,\partial_j^u,\partial_i^u)\\
 &-&2R_{2,s,F}(\partial_i^u,\nabla_{\partial_i^u}\partial_j^u,\partial_j^u,\partial_i^u)\\
&=&f_i^{\prime\prime\prime}+2u_i+2R_{2,s,F}(\textstyle\sum_{1\le k\le s}u_k\partial_k^t,\partial_j^u,\partial_j^u,\partial_i^u)+0
=f_i^{\prime\prime\prime}+4u_i\,.
\end{eqnarray*}
Assertion (3c) now follows.\end{proof}

\begin{proof}[Proof of Theorem \ref{thm-1.5} (1)-(3)] Adopt the notation of Section \ref{sect-1.3}. Let $1\le i\le r-1$ and
let $g=g_{3,r,\psi}$. We compute that the non-zero Christoffel symbols of the second kind are
$$\begin{array}{ll}
g(\nabla_{\partial_x}\partial_x,\partial_{u_r})=\psi^{\prime}(u_r),&
 g(\nabla_{\partial_x}\partial_{u_r},\partial_x)
=g(\nabla_{\partial_{u_r}}{\partial_x},\partial_x)=-\psi^{\prime}(u_r),\\
g(\nabla_{\partial_x}\partial_x,\partial_{u_i})=v_{i+1},&
 g(\nabla_{\partial_x}\partial_{u_i},\partial_x)=g(\nabla_{\partial_{u_i}}{\partial_x},\partial_x)=-v_{i+1},\\
g(\nabla_{\partial_x}\partial_x,\partial_{v_{i+1}})=u_i,&
 g(\nabla_{\partial_x}\partial_{v_{i+1}},\partial_x)=g(\nabla_{\partial_{v_{i+1}}}{\partial_x},\partial_x)=-u_i\,.
\end{array}$$
Consequently the non-zero Christoffel symbols of the first kind are
\begin{eqnarray}\nonumber
&&\nabla_{\partial_x}\partial_x=u_1\partial_{u_2}+...+u_{r-1}\partial_{u_r}
+v_2\partial_{v_1}+...+v_{r}\partial_{v_{r-1}}
+\psi^{\prime}(u_r)\partial_{v_r},\\
&&\nabla_{\partial_x}\partial_{u_i}=\nabla_{\partial_{u_i}}\partial_x=-v_{i+1}\partial_y,\quad\text{and}\quad
\nabla_{\partial_x}\partial_{v_{i+1}}=\nabla_{\partial_{v_{i+1}}}\partial_x=-u_i\partial_y\,.\label{eqn-2.b}
\end{eqnarray}
To apply Lemma \ref{lem-2.1}, we set
$$z_0=x,\ z_1=u_1,\ ...,\ z_r=u_r,\ z_{r+1}=v_r,\ ...,\ z_{2r}=v_1,\ z_{2r+1}=y\,.$$
Assertions (1) and (2) follow. We have
$$\begin{array}{ll}
\nabla_{\partial_{u_r}}\nabla_{\partial_x}\partial_x=\psi^{\prime\prime}\partial_{v_r},&
   \nabla_{\partial_x}\nabla_{\partial_{u_r}}\partial_x=0,\\
\nabla_{\partial_{u_i}}\nabla_{\partial_x}\partial_x=\partial_{u_{i+1}},&
   \nabla_{\partial_x}\nabla_{\partial_{u_i}}\partial_x=0,\\
\nabla_{\partial_{v_{i+1}}}\nabla_{\partial_x}\partial_x=\partial_{v_i},&
   \nabla_{\partial_x}\nabla_{\partial_{v_{i+1}}}\partial_x=0\,.
\end{array}$$
Assertions (3a) and (3b) follow. Assertion (3c) follows from these calculations and from Equation
(\ref{eqn-2.b}).
\end{proof}

\section{Model Spaces}\label{sect-3}

Throughout this section, we shall only list (possibly) non-zero entries of $g$, $R_g$, and of $\nabla R_g$ up to the usual
$\mathbb{Z}_2$ symmetries. We show $\mathcal{U}$ is a $0$-model for $\mathcal{M}$
by exhibiting a basis for $T_PM$ with the required normalizations for any $P\in M$.

\subsection{A $0$-model for $\mathcal{M}_{1,p,f}$}\label{sect-3.1}
 Choose a basis $\{X_1,...,X_p\}$ for
$\operatorname{Span}\{\partial_1^x,...,\partial_p^x\}$ in $T_PM$ so that
$H_f(X_i,X_j)=\delta_{ij}$. Expand $X_i=\sum_{1\le j\le p}\xi_{ij}\partial_j^x$ and let $\xi^{ij}$ be the
inverse matrix. Set
$Y_i:=\sum_{1\le j\le p}\xi^{ji}\partial_j^y$. Then
\begin{eqnarray*}
&&g_{1,p,f}(X_i,X_j)=c_{ij},\ \  g_{1,p,f}(X_i,Y_j)=\delta_{ij},\
\ \text{and}\\
&&R_{1,p,f}(X_i,X_j,X_k,X_l)=\delta_{il}\delta_{jk}-\delta_{ik}\delta_{jl},\
\end{eqnarray*}
where
$c_{ij}=c_{ji}$. Set
$\bar X_i=X_i-\ffrac12\textstyle\sum_{1\le j\le p}c_{ij}Y_j$ and $\bar Y_i:=Y_i$.
We may then conclude $\mathcal{U}_{1,p}$ is a $0$-model for $\mathcal{M}_{1,p,f}$ since
\begin{eqnarray*}
&&g_{1,p,f}(\bar X_i,\bar X_j)=0,\ \ g_{1,p,f}(\bar X_i,\bar
Y_j)=\delta_{ij},\ \ \text{and}\\
&&R_{1,p,f}(\bar X_i,\bar X_j,\bar X_k,\bar
X_l)=\delta_{il}\delta_{jk}-\delta_{ik}\delta_{jl}\,.
\end{eqnarray*}

\subsection{A $0$-model for $\mathcal{M}_{2,s,F}$}\label{sect-3.2}
Fix $P\in\mathbb{R}^{3s}$.
Define a new basis
for
$T_PM$ by setting:
\begin{equation}\label{eqn-3.a}
U_i:=\partial_i^u+\varepsilon_{i}\partial_i^t+\varrho_{i}\partial_i^v,\ \
 T_i:=\partial_i^t+\varepsilon_{i}\partial_i^v,\ \ \text{and}\ \ V_i:=\partial_i^v\end{equation}
where the constants $\varepsilon_i$ and $\varrho_i$ will be specified below.
Let $i\ne j$. Then:
\begin{eqnarray*}
&&g_{2,s,F}(U_i,T_i)=\varepsilon_{i}-\varepsilon_{i}=0,\ \
g_{2,s,F}(U_i,U_i)=g_{2,s,F}(\partial_i^u,\partial_i^u)-\varepsilon_{i}^2+2\varrho_{i},\\
&&g_{2,s,F}(T_i,T_i)=-1,\ \ g_{2,s,F}(U_i,V_i)=1,\ \
R_{2,s,F}(U_i,U_j,U_j,T_i)=1,\ \ \text{and}\\
&&R_{2,s,F}(U_i,U_j,U_j,U_i)=(\partial_i^u)^2f_i+(\partial_j^u)^2f_j+
|u|^2+2\varepsilon_i+2\varepsilon_j\,.
\end{eqnarray*}
We set
$$\varepsilon_{i}:=-\ffrac12(\partial_i^u)^2f_i-\ffrac14|u|^2
 \ \ \text{and}\ \ \varrho_i:=\ffrac12\{\varepsilon_{i}^2-g_{2,s,F}(\partial_i^u,\partial_i^u)\}\,.
$$
As $g_{2,s,F}(U_i,U_i)=R_{2,s,F}(U_i,U_j,U_j,U_i)=0$, $\mathcal{U}_{2,s}$
is a $0$-model for $\mathcal{M}_{2,s,F}$.

\subsection{A $0$-model for $\mathcal{M}_{3,r,\psi}$}\label{sect-3.3}
Let $\varepsilon_i$ be real parameters to be specified below.
Define a new basis $\{X,Y,U_1,...,U_r,V_1,...,V_r\}$ for
$T_P\mathbb{R}^{2r+2}$ by setting:
 $$X=\varepsilon_0\{\partial_x-\ffrac12g_{3,r,\psi}(\partial_x,\partial_x)\partial_y\},\
Y=\varepsilon_0^{-1}\partial_y,\ U_i=\varepsilon_i\partial_{u_i},\
V_i=\varepsilon_i^{-1}\partial_{v_i}\,. $$ The non-zero entries in
$g_{3,r,\psi}$ are given by $g_{3,r,\psi}(X,Y)=1$ and
$g_{3,r,\psi}(U_i,V_i)=1$. We apply Theorem \ref{thm-1.5} (3) to
see the non-zero entries in $R_{3,r,\psi}$ and $\nabla
R_{3,r,\psi}$ are
\begin{eqnarray*}
&&R_{3,r,\psi}(X,U_r,U_r,X)=\varepsilon_0^2\varepsilon_r^2\psi^{\prime\prime}(u_r),\\
&&R_{3,r,\psi}(X,U_i,V_{i+1},X)=\varepsilon_0^2\varepsilon_i\varepsilon_{i+1}^{-1}\ \ \text{for}\ \ 1\le
i\le r-1,\\
&&\nabla R_{3,r,\psi}(X,U_r,U_r,X;U_r)=\varepsilon_0^2\varepsilon_r^3\psi^{\prime\prime\prime}(u_r)\,.
\end{eqnarray*}
Assume $\psi^{\prime\prime}>0$. We can show that
$\mathcal{U}_{3,r}$ is a $0$-model for $\mathcal{M}_{3,r,\psi}$ by
setting:
 $$\varepsilon_r=(\psi^{\prime\prime})^{-1/2},\ \
\varepsilon_0=1,\ \ \text{and}\ \
 \varepsilon_i=\varepsilon_r\ \ \text{for}\ \ 1\le i<r\,.$$
If in addition we suppose that $\psi^{\prime\prime\prime}\ne0$,
then more is true. We show $\mathcal{U}_{3,r}^1$is a $1$-model for
$\mathcal{M}_{3,r,\psi}$ by setting:
$$\varepsilon_r=\psi^{\prime\prime}(\psi^{\prime\prime\prime})^{-1},\
\
 \varepsilon_0=(\varepsilon_r^2\psi^{\prime\prime})^{-1/2},\ \
 \text{and}\ \
 \varepsilon_i=\varepsilon_0^{-2}\varepsilon_{i+1}\ \ \text{for}\ \ 1\le i<r\,.$$

\section{Irreducibility}\label{sect-4}
\subsection{The model $\mathcal{M}_{1,p,f}$} We adopt the notation of Section \ref{sect-1.2}. Let
$B_{1,p}$ be the algebraic curvature tensor on
$\mathbb{R}^p=\operatorname{Span}\{X_1,...,X_p\}$ defined by
$$
B_{1,p}(X_i,X_j,X_k,X_l)=\delta_{il}\delta_{jk}-\delta_{ik}\delta_{jl}\,.$$

\begin{lemma}\label{lem-4.1}\
\begin{enumerate}
\item Let $0\ne \xi_1\in\mathbb{R}^p$. If $B_{1,p}(\xi_1,\xi_2,\eta_1,\eta_2)=0$ $\forall$
$\eta_1,\eta_2\in\mathbb{R}^p$, then
$\xi_2=\lambda \xi_1$.
\item $(\mathbb{R}^p,B_{1,p})$ is irreducible.
\end{enumerate}\end{lemma}

\begin{proof} Let
$g_0$ be the usual Euclidean inner product;
$g_0(X_i,X_j):=\delta_{ij}$. Then:
$$B_{1,p}(\eta_1,\eta_2,\eta_3,\eta_4)=g_0(\eta_1,\eta_4)g_0(\eta_2,\eta_3)-g_0(\eta_1,\eta_3)g_0(\eta_2,\eta_4)\,.$$
Let $O(p)$ be the usual Euclidean orthogonal group. If $\theta\in
O(p)$, then $\theta^*B_{1,p}=B_{1,p}$. By applying a suitable
element of $\theta\in O(p)$ and rescaling if necessary, we may
assume without loss of generality $\xi_1=X_1$ in proving Assertion
(1). We expand $\xi_2=\sum_{1\le i\le p}a_iX_i$. For $i>1$,
$a_i=R(\xi_1,\xi_2,X_i,X_1)=0$. Assertion (1) follows.

Suppose that we have a non-trivial decomposition $\mathbb{R}^p=W_1\oplus W_2$ which induces a
decomposition $B_{1,p}=B_{1,p}^1\oplus B_{1,p}^2$. Let $0\ne \xi_i\in W_i$. Then
$B_{1,p}(\xi_1,\xi_2,\cdot,\cdot)=0$ so, by Assertion (1), $\xi_1=\lambda \xi_2$; this is false.
\end{proof}

\begin{proof}[Proof of Theorem \ref{thm-1.3} (4)] We showed in Section \ref{sect-3.1} that $\mathcal{U}_{1,p}$ is a $0$-model
for
$\mathcal{M}_{1,p,f}$. Thus we must only show that
$\mathcal{U}_{1,p}$ is irreducible. Let
\begin{eqnarray*}
&&K:=\operatorname{Span}\{Y_1,...,Y_p\}
=\{\eta\in\mathbb{R}^{2p}:R(\xi_1,\xi_2,\xi_3,\eta)=0\ \ \forall\ \
\xi_i\in\mathbb{R}^{2p}\}\,.
\end{eqnarray*}
Let $\pi$ be the natural projection from $\mathbb{R}^{2p}$ to
$\mathbb{R}^p=\operatorname{Span}\{X_1,...,X_p\}=\mathbb{R}^{2p}/K$.
We then have that
$A_{1,p}=\pi^*B_{1,p}$.
Suppose there is a non-trivial decomposition:
\begin{equation}\label{eqn-4.a}
\mathbb{R}^{2p}=V_1\oplus V_2,\ \ g_{1,p}=g_{1,p}^1\oplus
g_{1,p}^2,\ \ \text{and}\ \ A_{1,p}=A_{1,p}^1\oplus A_{1,p}^2\,.
\end{equation}

We argue for a contradiction. Since $V_1\perp V_2$, the metrics
$g_{1,p}^i$ are non-trivial on $V_i$. In particular, the subspaces
$V_i$ are not totally isotropic. Equation (\ref{eqn-4.a}) induces
a corresponding decomposition $$\mathbb{R}^p=V_1/\{K\cap
V_1\}\oplus V_2/\{K\cap V_2\}\ \ \text{and}\ \
B_{1,p}=B_{1,p}^1\oplus B_{1,p}^2\,. $$ By Assertion (1), this
decomposition of $\mathbb{R}^p$ is trivial; we assume that the
notation is chosen so $V_2/\{K\cap V_2\}=\{0\}$ and hence
$V_2\subset K$ so $V_2$ is totally isotropic. This is a
contradiction.\end{proof}

\subsection{The model $\mathcal{M}_{2,s,F}$} We adopt the notation of Section \ref{sect-1.2}. We define an algebraic
curvature tensor $B_{2,s}$ on
$\mathbb{R}^{2s}:=\operatorname{Span}\{U_1,...,U_s,T_1,...,T_s\}$
by setting:
$$B_{2,s}(U_i,U_j,U_k,T_l)=\delta_{il}\delta_{jk}-\delta_{ik}\delta_{jl}\,.
$$
\begin{lemma}\label{lem-4.2}\
\begin{enumerate}
\item  Let $0\ne \xi_1,\xi_2\in\mathbb{R}^{2s}$. If $B_{2,s}(\xi_1,\xi_2,\eta_1,\eta_2)=0$ and
$B_{2,s}(\xi_1,\eta_1,\eta_2,\xi_2)=0$ for all
$\eta_1,\eta_2\in\mathbb{R}^{2s}$, then $\xi_1,\xi_2\in\operatorname{Span}\{T_1,...,T_s\}$.
\item  $(\mathbb{R}^{2s},B_{2,s})$ is irreducible.
\end{enumerate}\end{lemma}

\begin{proof} We extend
$\theta\in O(s)$ to act diagonally on
$\mathbb{R}^{2s}=\mathbb{R}^s\oplus\mathbb{R}^s$; we then have $\theta^*B_{2,s}=B_{2,s}$. Suppose that
$\xi_1\not\in\operatorname{Span}\{T_1,...,T_s\}$. By applying a suitably chosen element $\xi\in
O(s)$ and rescaling if necessary, we may assume without loss of generality
\begin{eqnarray*}
&&\xi_1=U_1+b_1T_1+...+b_sT_s\ \ \text{and}\\
&&\xi_2=c_1U_1+...+c_sU_s+d_1T_1+...+d_sT_s
\end{eqnarray*}
for suitably chosen constants $\{b_1,...,b_s,c_1,...,c_s,d_1,...,d_s\}$. Let $i>1$. We have
$$0=B_{2,s}(\xi_1,\xi_2,U_i,T_1)=c_i\ \ \text{and}\ \
0=B_{2,s}(\xi_1,U_i,T_i,\xi_2)=c_1\,.$$
This shows that $c_1=0$ and $c_i=0$ so $\xi_2=d_1T_1+...+d_sT_s$. Furthermore,
$$0=B_{2,s}(\xi_1,\xi_2,U_i,U_1)=d_i\ \ \text{and}\ \ 0=B_{2,s}(\xi_1,U_i,U_i,\xi_2)=d_1\,.$$
This implies $\xi_2=0$ which is a contradiction. Thus $\xi_1\in\operatorname{Span}\{T_1,...,T_s\}$. As the roles of
$\xi_1$ and $\xi_2$ are symmetric, we conclude $\xi_2\in\operatorname{Span}\{T_1,...,T_s\}$ as well;
Assertion (1) follows.

Suppose given a non-trivial decomposition $\mathbb{R}^{2s}=W_1\oplus W_2$ which induces a
decomposition $B_{2,s}=B_{2,s}^1\oplus B_{2,s}^2$. Choose $0\ne \xi_i\in W_i$. By Assertion (2),
$\xi_1,\xi_2$ belong to $\operatorname{Span}\{T_i\}$. Thus
$W_1\subset\operatorname{Span}\{T_i\}$ and $W_2\subset\operatorname{Span}\{T_i\}$. Thus
$\mathbb{R}^{2s}$  is contained in $\operatorname{Span}\{T_i\}$ which is false.
\end{proof}

\begin{proof}[Proof of Theorem \ref{thm-1.4} (4)] We showed in Section \ref{sect-3.2} that $\mathcal{U}_{2,s}$ is a $0$-model
for
$\mathcal{M}_{2,s,F}$. Thus it suffices to show that
$\mathcal{U}_{2,s}$ is irreducible. Let
\begin{eqnarray*}
&&L:=\operatorname{Span}\{V_1,...,V_s\}
=\{\eta\in\mathbb{R}^{3s}:A_{2,s}(\xi_1,\xi_2,\xi_3,\eta)=0\ \ \forall\ \
\xi_i\in\mathbb{R}^{3s}\}\,.
\end{eqnarray*}
Let $\pi$ be the natural projection from $\mathbb{R}^{3s}$ to
$$\mathbb{R}^{2s}:=\operatorname{Span}\{U_1,...,U_p,T_1,...,T_p\}=\mathbb{R}^{3s}/L\,.$$ We have
$A_{2,s}=\pi^*B_{2,s}$.
Suppose we have a non-trvial decomposition
\begin{equation}\label{eqn-4.b}
\mathbb{R}^{3s}=V_1\oplus V_2,\ \ g_{2,s}=g_{2,s}^1\oplus g_{2,s}^2,\ \ \text{and}\ \
A_{2,s}=A_{2,s}^1\oplus A_{2,s}^2\,.
\end{equation}

We argue for a contradiction. We argue as above to see $V_1$ and $V_2$ are not totally isotropic. Equation (\ref{eqn-4.b})
induces a corresponding decomposition
$$\mathbb{R}^{2s}=V_1/\{L\cap V_1\}\oplus V_2/\{L\cap V_2\}\ \ \text{and}\ \ B_{2,s}=B_{2,s}^1\oplus
B_{2,s}^2\,.$$ By Lemma \ref{lem-4.2}, this decomposition must be trivial. We assume the notation chosen so
that
$V_2\subset L$. Thus $V_2$ is totally isotropic. This is a contradiction.\end{proof}

\subsection{The model $\mathcal{U}_{3,r}$} Adopt the notation of Section \ref{sect-1.3}. Let $1\le i<r$. The non-zero
entries in the curvature operator are, up to the usual $\mathbb{Z}_2$ symmetries,
\begin{equation}\label{eqn-4.c}
\begin{array}{ll}
A_{3,r}(X,U_r)U_r=Y,& A_{3,r}(X,U_r)X=-V_r,\\
A_{3,r}(X,U_i)V_{i+1}=Y,& A_{3,r}(X,U_i)X=-U_{i+1},\\
A_{3,r}(X,V_{i+1})U_i=Y,& A_{3,r}(X,V_{i+1})X=-V_i\,.
\end{array}\end{equation}
If $\xi\in\mathbb{R}^{2r+2}$, then:
\begin{eqnarray*}
&&J(\xi)X\in\operatorname{Span}\{U_2,...,U_r,V_1,...,V_r,Y\},\\
&&J(\xi)U_r\in\operatorname{Span}\{V_r,Y\},\qquad
J(\xi)U_i\in\operatorname{Span}\{U_{i+1},Y\},\\
&&J(\xi)V_i\in\operatorname{Span}\{V_{i-1},Y\},\phantom{a...}
J(\xi)Y=J(\xi)V_1=0\,.
\end{eqnarray*}
\begin{proof}[Proof of Theorem \ref{thm-1.6} (3)] Display (\ref{eqn-4.c}) shows
$J(\xi)^{2r}=0$. As $J(X)^{2r-1}U_1=V_1$, $\mathcal{U}_{3,r}$ is
$2r$-Osserman nilpotent; Theorem \ref{thm-1.6} (3)
follows.\end{proof}

\begin{proof}[Proof of Theorem \ref{thm-1.5} (4)] We showed in Section \ref{sect-3.3} that $\mathcal{U}_{3,r}$ is a $0$-model
for
$\mathcal{M}_{3,r,\psi}$. Thus it suffices to show $\mathcal{U}_{3,r}$ is irreducible. We suppose the contrary and argue for a
contradiction. Suppose there is a non-trivial decomposition
\begin{equation}\label{eqn-4.d}
\mathbb{R}^{2r+2}=W_1\oplus W_2,\quad g_{3,r}=g_{3,r}^1\oplus
g_{3,r}^2,
 \quad\text{and}\quad A_{3,r}=A_{3,r}^1\oplus A_{3,r}^2\,.
\end{equation}
As above, neither $W_1$ nor $W_2$ can be totally isotropic.
Decompose $X=X_1+X_2$. Then either $J(X_1)$ or
$J(X_2)$ is nilpotent of order $2r$; we may assume without loss of generality that the
notation is chosen so that $J(X_1)$ is nilpotent of order $2r$. Since $J(X_1)X_1=0$,
this implies $\dim(W_1)\ge 2r+1$. Since the decomposition is assumed non-trivial, this implies
$\dim(W_2)=1$. Let $\xi$ span $W_2$; $\xi$ can not be a null vector since $W_2$ is not totally isotropic.
On the other hand since $\dim(W_2)=1$, $A_{3,r}(\eta_1,\eta_2)\xi=0$ for $\eta_i\in W_2$. The decomposition of
Equation (\ref{eqn-4.d}) then shows $A_{3,r}(\eta_1,\eta_2)\xi=0$ for all $\eta_1,\eta_2\in\mathbb{R}^{2r+2}$. This
implies
$\xi\in\operatorname{Span}\{V_1,Y\}$ which is a totally isotropic subspace; this is a contradiction.\end{proof}

\section{Homogeneity}\label{sect-5}
\subsection{The manifolds $\mathcal{M}_{1,p,f}$} If $\phi$ is a symmetric bilinear form on $V$, then we may
define an algebraic curvature tensor $R(\phi)$ on $V$ by setting:
\begin{equation}\label{eqn-5.a}
R(\phi)(\xi_1,\xi_2,\xi_3,\xi_4):=
 \phi(\xi_1,\xi_4)\phi(\xi_2,\xi_3)-\phi(\xi_1,\xi_3)\phi(\xi_2,\xi_4)\,.
\end{equation}
One then has, see for example the discussion in \cite{refDG},
\begin{lemma}\label{lem-5.1}
Let $\phi_1$ and $\phi_2$ be symmetric positive definite bilinear forms on a vector
space $V$ of dimension at least $3$. If
$R({\phi_1})=R({\phi_2})$, then $\phi_1=\phi_2$.
\end{lemma}

\begin{proof}[Proof of Theorem \ref{thm-1.3} (6)] Adopt the notation of Section \ref{sect-1.2}. Fix $P\in\mathbb{R}^{2p}$ and
let
$V_P:=T_P\mathbb{R}^{2p}$. We consider a
$1$-model
$$\mathcal{V}_P^1:=(V_P,g_{1,p,f}|_{V_P},R_{1,p,f}|_{V_P},\nabla R_{1,p,f}|_{V_P})\,.$$
Also consider the subspace
$$Y_P:=\{\eta\in V_P:R_{1,p,f}(\xi_1,\xi_2,\xi_3,\eta)=0\ \ \forall\ \
\xi_i\in V_P\}
=\operatorname{Span}\{\partial_1^y,...,\partial_p^y\}\,.
$$
Let $\pi$ be the natural projection from $V_P$ to $X_P:=V_P/Y_P$. As
$$H(\xi_1,\xi_2)=0,\ \ R_{1,p,f}(\xi_1,\xi_2,\xi_3,\xi_4)=0,\ \ \text{and}\ \
\nabla R_{1,p,f}(\xi_1,\xi_2,\xi_3,\xi_4;\xi_5)=0$$
if any $\xi_i\in Y_P$, there are structures
$H_{X,P}$,
$A_{X,P}$, and
$\AA_{X,P}$ on $X_P$ which are characterized by the identities:
$$\pi^*H_{X,P}=H_f|_{V_P},\quad\pi^*A_{X,P}=R_{1,p,f}|_{V_P},\quad\text{and}\quad
 \pi^*\AA_{X,P}=\nabla R_{1,p,f}|_{V_P}\,.
$$

Assume $\mathcal{M}_{1,p,f}$ is $1$-curvature homogeneous. Let $P,Q\in\mathbb{R}^{2p}$. Let $\Theta$ be an isomorphism
from
$\mathcal{V}_P^1$ to
$\mathcal{V}_Q^1$. It is immediate from the defining relation that
$\Theta(Y_P)\subset Y_Q$; a dimension count then implies $\Theta(Y_P)=Y_Q$. Consequently
$\Theta$ induces a map $\tilde\Theta:X_P\rightarrow X_Q$ so
$$\tilde\Theta^*A_{X,Q}=A_{X,P}\quad\text{and}\quad\tilde\Theta^*\AA_{X,Q}=\AA_{X,P}\,.$$
We adopt the notation of Equation (\ref{eqn-5.a}) and let $R(\phi)$ be the curvature tensor defined by a bilinear form $\phi$.
Since
$$R(H_{X,P})=A_{X,P}=\tilde\Theta^*(A_{X,Q})=R(\tilde\Theta^*H_{X,Q}),$$
Lemma \ref{lem-5.1} implies that
$H_{X,P}=\tilde\Theta^*H_{X,Q}$. Let $||_\phi^2$ denote the norm taken with respect to a positive definite bilinear form
$\phi$. We then have
$$\alpha_1(P)=||A^1_{X,P}||_{H_{X,P}}^2=||A^1_{X,Q}||_{H_{X,Q}}^2=\alpha_1(Q)\,.$$
Consequently $\alpha_1$ is
constant. \end{proof}

\subsection{The manifolds $\mathcal{M}_{2,s,F}$} We begin by studying the Lie group associated to the model
$\mathcal{U}_{2,s}$. We say that
$\BB=\{u_1,...,u_s,t_1,...,t_s,v_1,...,v_s\}$ is a {\it normalized} basis for
$\mathbb{R}^{3s}$ if the normalizations of Equation (\ref{eqn-1.a}) hold, i.e.
\begin{eqnarray*}
&&g_{2,s}(u_i,v_j)=\delta_{ij},\ g_{2,s}(t_i,t_j)=-\delta_{ij},\ \ \text{and}\\
&&A_{2,s}(u_i,u_j,u_k,t_l)=\delta_{il}\delta_{jk}-\delta_{ik}\delta_{jl}\,.
\end{eqnarray*}
Let $O(s)\subset M_s(\mathbb{R})$ be the usual orthogonal group of $s\times s$ matrices;
$\kappa_{ij}\in O(s)$ if and only if $\sum_k\kappa_{ik}\kappa_{jk}=\delta_{ij}$.

\begin{lemma}\label{lem-5.2}
 Let $\BB$ and $\tilde\BB$ be two normalized bases for $\mathbb{R}^{3s}$. If $s\ge3$, then there exists a matrix $\kappa_1\in
O(s)$ and matrices $\kappa_2,\ \kappa_3,\ \kappa_5\in M_s(\mathbb{R})$ so that:
\begin{eqnarray*}
&&\tilde u_i=\textstyle\sum_j\{\kappa_{1,ij}u_j+\kappa_{2,ij}t_j+\kappa_{3,ij}v_j\},\\
&&\tilde t_i=\textstyle\sum_j\{\kappa_{1,ij}t_j+\kappa_{5,ij}v_j\},\quad\text{and}\quad
  \tilde v_i=\textstyle\sum_j\kappa_{1,ij}v_j\,.
\end{eqnarray*}
\end{lemma}

\begin{proof} We note that
\begin{eqnarray*}
Y:&=&\{\eta\in\mathbb{R}^{3s}:R_{2,s}(\zeta_1,\zeta_2,\zeta_3,\eta)=0\ \text{for all}\
\zeta_1,\zeta_2,\zeta_3
\in\mathbb{R}^3\}\\
&=&\operatorname{Span}\{v_1,...,v_s\}=\operatorname{Span}\{\tilde v_1,...,\tilde v_s\},\quad\text{and}\\
Y^\perp:&=&\{\eta\in\mathbb{R}^{3s}:g_{2,s}(\eta,\zeta)=0\ \text{for all}\ \zeta\in Y\}\\
&=&\operatorname{Span}\{t_1,...,t_s,v_1,...,v_s\}=\operatorname{Span}\{\tilde t_1,...,\tilde t_s,\tilde v_1,...,\tilde v_s\}\,,
\end{eqnarray*}
Consequently we may express
\begin{eqnarray*}
&&\tilde u_i=\textstyle\sum_j\{\kappa_{1,ij}u_j+\kappa_{2,ij}t_j+\kappa_{3,ij}v_j\},\\
&&\tilde t_i=\textstyle\sum_j\{\kappa_{4,ij}t_j+\kappa_{5,ij}v_j\},\quad\text{and}\quad
\tilde v_i=\sum_j\kappa_{6,ij}v_j\,.
\end{eqnarray*}
We verify that $\kappa_4\in O(s)$ by checking
$$-\delta_{ij}=g_{2,s}(\tilde t_i,\tilde t_j)=\textstyle\sum_{k,l}\kappa_{4,ik}\kappa_{4,jl}g_{2,s}(t_i,t_j)
=-\textstyle\sum_k\kappa_{4,ik}\kappa_{4,jk}\,.$$

The orthogonal group acts diagonally on $\mathbb{R}^{3s}$ by
$$\kappa:u_i\rightarrow\textstyle\sum_j\kappa_{ij}u_j,\ \
  \kappa:t_i\rightarrow\textstyle\sum_j\kappa_{ij}t_j,\ \ \text{and}\ \
  \kappa:v_i\rightarrow\textstyle\sum_j\kappa_{ij}v_j\,.
$$
This action preserves the structures involved. By making a suitable change of basis, therefore, we may suppose without loss of
generality that $\kappa_4=\operatorname{id}$ in the proof of the Lemma, i.e. that we have:
\begin{eqnarray*}
&&\tilde u_i=\textstyle\sum_j\{\kappa_{1,ij}u_j+\kappa_{2,ij}t_j+\kappa_{3,ij}v_j\},\\
&&\tilde t_i=t_i+\textstyle\sum_j\kappa_{5,ij}v_j,\quad \text{and}\quad \tilde v_i=\sum_j\kappa_{6,ij}v_j\,.
\end{eqnarray*}

To show $\kappa_1=\operatorname{id}$, fix $i\ne j$. Since $s\ge3$,
we can choose $$0\ne u\in\operatorname{Span}_{k\ne
j}\{u_k\}\cap\operatorname{Span}_{k\ne j;1\le\ell\le s}\{\tilde
u_k,\tilde t_\ell,\tilde v_\ell\}\,.$$ Expand $u=\sum_{k\ne
j}\varepsilon_ku_k$. As $\tilde t_j=t_j+\sum_k\kappa_{5,jk}v_k$,
we may compute
\begin{eqnarray*}
0&=&R_{2,s}(\tilde u_i,u,u,\tilde
t_j)=R_{2,s}(\textstyle\sum_k\kappa_{1,ik}u_k,\sum_{a\ne
j}\varepsilon_au_a,\sum_{b\ne j}\varepsilon_bu_b,t_j)\\
&=&\kappa_{1,ij}\textstyle\sum_{a\ne j}\varepsilon_a^2\,.
\end{eqnarray*}
This shows $\kappa_{1,ik}=0$ for $i\ne k$ so $\kappa_1$ is diagonal. Since
$$1=R_{2,s}(\tilde u_i,\tilde u_j,\tilde u_j,\tilde t_i)=\kappa_{1,ii}\kappa_{1,jj}\kappa_{1,jj},$$
and similarly $1=\kappa_{1,jj}\kappa_{1,ii}\kappa_{1,ii}$, we have $\kappa_{1,ii}=1$ as desired. The identity
$g_{2s}(\tilde u_i,\tilde v_j)=\delta_{ij}$ then shows $\kappa_6=\operatorname{id}$ in this special situation.
\end{proof}

Fix $P\in\mathbb{R}^{3s}$ and let $V_P:=T_P\mathbb{R}^{3s}$.
Consider a $1$-model
$$\mathcal{V}_P^1:=(V_P,g_{2,s,F}|_{V_P},R_{2,s,F}|_{V_P},\nabla
R_{2,s,F}|_{V_P})\,. $$ Also consider the subspaces
\begin{eqnarray*}
&&Y_P:=\{\eta\in V_P:R_{2,s,F}(\xi_1,\xi_2,\xi_3,\eta)=0\ \forall\ \xi_i\in V_P\}
=\operatorname{Span}\{\partial_1^v,...,\partial_s^v\},\text{ and}\\
&&Y_P^\perp=\{\eta\in\mathbb{R}^{3s}:g_{2,s,F}(\eta,\xi_1)=0\ \forall\ \xi_1\in Y_P\}
=\operatorname{Span}\{\partial_1^t,...,\partial_s^t,\partial_1^v,...,\partial_s^v\}\,.
\end{eqnarray*}
Let $\pi$ be the natural projection from $V_P$ to  $X_P:=V_P/Y_P$. There is a natural covariant derivative algebraic curvature
tensor $A_{X,P}^1$ on $X_P$ so
$\pi^*A_{X,P}^1=\nabla R_{2,s,F}|_{V_P}$;
$$A_{X,P}^1(\tilde U_i,\tilde U_j,\tilde U_j,\tilde U_i;\tilde U_i):=(\partial_i^u)^3f_i+4u_i\,.$$

The elements $\{\tilde U_i:=\pi
\partial_i^u\}$ are a basis for $X_P$.
Define a non-degenerate bilinear form $L_P$ on
$X_P$ by requiring that
$$L_P(\tilde U_i,\tilde U_j)=\delta_{ij}\,.$$
If $\Theta$ is an isomorphism from $\mathcal{V}_P^1$ to $\mathcal{V}_Q^1$, then clearly $\Theta(Y_P)=Y_Q$. Consequently
$\Theta(Y_P^\perp)=Y_Q^\perp$ so
$\Theta$ induces a map
$\tilde\Theta$ from $X_P$ to $X_Q$.

To construct the normalized basis of Equation (\ref{eqn-3.a}) we set:
$$U_i=\partial_i^u+\operatorname{Span}_j\{\partial_j^t,\partial_j^v\},\
  T_i=\partial_i^t+\operatorname{Span}_j\{\partial_j^v\},\text{ and }V_i=\partial_i^v\,.
$$
We apply Lemma \ref{lem-5.2} to expand
$$\Theta
U_i=\textstyle\sum_j\kappa_{1,ij}U_i
+\operatorname{Span}_j\{T_j,V_j\}$$
where
$\kappa_1\in O(s)$. Since $\tilde U_i=\pi U_i$,  $\Theta\tilde U_i=\sum_j\kappa_{1,ij}\tilde U_j$. The following Lemma
is now immediate:
\begin{lemma}\label{lem-5.3}
If $\Theta$ is an isomorphism from $\mathcal{V}_P^1$ to $\mathcal{V}_Q^1$, then
$\tilde\Theta^*L_Q=L_P$.\end{lemma}

\begin{proof}[Proof of Theorem \ref{thm-1.4} (6)]
Assume $\mathcal{M}_{2,s,F}$ is $1$-curvature homogeneous. Let $P$ and $Q$ be points of $\mathbb{R}^{3s}$. Let $\Theta$ be an
isomorphism from
$\mathcal{V}_P^1$ to $\mathcal{V}_Q^1$.  Since
$\alpha_2=|A_{X,P}^1|_{L_P}^2$, Lemma \ref{lem-5.3} implies $\alpha_2$ must be constant.\end{proof}

\subsection{The manifolds $\mathcal{M}_{3,r,\psi}$} Adopt the notation of Section \ref{sect-1.4}. Assume that
$\psi^{\prime\prime}>0$ and that
$\psi^{\prime\prime\prime}>0$. Set
$$K_P:=\{\xi\in\mathbb{R}^{2r+2}:\exists\xi_i\in T_P\mathbb{R}^{2r+2}
\text{ so }\nabla^2R_{3,r,\psi}(\xi_1,\xi_2,\xi_3,\xi_4;\xi_5,\xi)\ne0\}\,.
$$

\begin{proof}[Proof of Theorem \ref{thm-1.5} (6)] Assume that $\psi^{\prime\prime}>0$ and $\psi^{\prime\prime\prime}>0$ for
all points of $\mathbb{R}$. The possibly non-zero entries in
$\nabla^2R_{3,r,\psi}$ are given by:
\begin{eqnarray*}
&&\nabla^2R_{3,r,\psi}(\partial_x,\partial_{u_r},\partial_{u_r},\partial_x;\partial_x,\partial_x)
 =u_{r-1}\psi^{\prime\prime\prime}(u_r),\\
&&\nabla^2R_{3,r,\psi}(\partial_x,\partial_{u_r},\partial_{u_r},\partial_x;\partial_{u_r},\partial_{u_r})
=\psi^{\prime\prime\prime\prime}(u_r)\,.
\end{eqnarray*}
We expand $\xi=\xi_0\partial_x+\xi_1\partial_1^u+...+\xi_r\partial_r^u+\xi_{r+1}\partial_r^v+...+\xi_{2r}\partial_1^v+
\xi_{2r+1}\partial_y$. Then
$$
K_P=\left\{\begin{array}{ll}
\{\xi\in\mathbb{R}^{2r+2}:\xi_0^2+\xi_r^2\ne0\}&\text{if }\psi^{\prime\prime\prime\prime}\ne0\text{ and }u_{r-1}\ne0,\\
\{\xi\in\mathbb{R}^{2r+2}:\xi_r\ne0\}&\text{if }\psi^{\prime\prime\prime\prime}\ne0\text{ and }u_{r-1}=0,\\
\{\xi\in\mathbb{R}^{2r+2}:\xi_0\ne0\}&\text{if }\psi^{\prime\prime\prime\prime}=0\text{ and }u_{r-1}\ne0,\\
\{0\}&\text{if }\psi^{\prime\prime\prime\prime}=0\text{ and }u_{r-1}=0\,.\\
\end{array}\right.$$

Suppose that $\mathcal{M}_{3,r,\psi}$ is curvature $2$-homogeneous. Then $K_P$ is diffeomorphic to $K_Q$ for any two points
$P$ and $Q$ of $\mathbb{R}^{3s}$. Let $P=(0,...,1,u_r,0,....,0)$ and
$Q=(0,...,0,u_r,0,...,0)$. Suppose $\psi^{\prime\prime\prime\prime}(u_r)\ne0$. Then $K_P$ is connected and $K_Q$ is not
connected; this is a contradiction. Suppose $\psi^{\prime\prime\prime\prime}(u_r)=0$. Then $K_P$ is non-empty and $K_Q$ is
empty; again, this is a contradiction. \end{proof}

\section{Symmetric Spaces}\label{sect-6}

\begin{proof}[Proof of Theorem \ref{thm-1.8}]
We extend an argument of E. Cartan's from the Riemannian
to the pseudo-Riemannian setting. Let
$\{e_i\}$ be a parallel frame field along a geodesic $\sigma$. Then
$$\partial_tR_{ijkl}(t)=\nabla R(e_i,e_j,e_k,e_l;\dot\sigma)=0\,.$$
Thus $R(e_i,\dot\sigma)\dot\sigma=c_{ij}e_j$ for suitably chosen
constants $c_{ij}$. Let
$Y(t)$ be a Jacobi vector field. Express $Y(t)=a_i(t)e_i(t)$. Then:
\begin{eqnarray*}
&&0=\ddot Y(t)+R(Y,\dot\sigma)\dot\sigma=\{\ddot
a_j(t)+\textstyle\sum_ja_i(t)c_{ij}\}e_j(t)\ \ \text{so}\\
&&0=\ddot a_j(t)+a_i(t)c_{ij}\ \ \text{for}\ \ 1\le j\le m\,.
\end{eqnarray*}
Since $-a_j(\xi;-t)$ still satisfies the Jacobi
equation with the same initial condition, $a_j(\xi;t)=-a_j(\xi;-t)$ so $a_j$
is an odd function of
$t$. Let $g_{ij}:=g(e_i,e_j)$ be independent of $t$. Then
$$g(Y_\xi(t),Y_\eta(t))=g_{ij}a_i(\xi;t)a_j(\eta;t)=g(Y_\xi(-t),Y_\eta(-t))$$
is an even
function of
$t$. Since the geodesic involution takes $Y_\xi(t)$ to $-Y_\xi(-t)$, this shows the
geodesic involution is an isometry and establishes the first assertion.

Let $P,Q\in M$. We suppose $P\ne Q$. Since $\exp_P$ is a diffeomorphism from $T_PM$ to
$M$, we can choose a geodesic $\sigma$ so $\sigma(0)=P$ and $\sigma(1)=Q$. Let
$R=\sigma(\frac12)$. Then the geodesic involution centered at $R$ interchanges $P$ and
$Q$ and is an isometry. Thus $M$ is a homogeneous space.
\end{proof}

\section*{Acknowledgments} Research of P. Gilkey partially supported by the
MPI (Leipzig). Research of S. Nik\v cevi\'c partially supported by MM 1646 (Serbia).

\end{document}